\documentclass[letterpaper,11pt]{article}
\usepackage{amsmath, amsthm, amssymb}    
\usepackage{bridges}			
\usepackage{graphicx}			
\usepackage[colorlinks=true, urlcolor=blue, citecolor=black, linkcolor=black]{hyperref}  
\usepackage{subcaption}			

\title{A Woven Klein Quartic}

\author{Chaim Goodman-Strauss
\\
National Museum of Mathematics; chaimgoodmanstrauss@gmail.com\\}

\begin{document}

\maketitle

\thispagestyle{empty}

\begin{abstract}

We describe a new method of weaving a model of the Klein quartic, a highly symmetric, but abstract genus-3 surface akin to a platonic polyhedron, with negatively-curved geometry. The Klein quartic cannot be realized in its fully symmetric form in three-dimensional space, but this model exhibits the most rigid symmetry that is possible. With remarkably little time and material you can have a Klein quartic of your own!
\end{abstract}

\begin{figure}[h!tbp]
	\centering
	\includegraphics[width=2.4in]{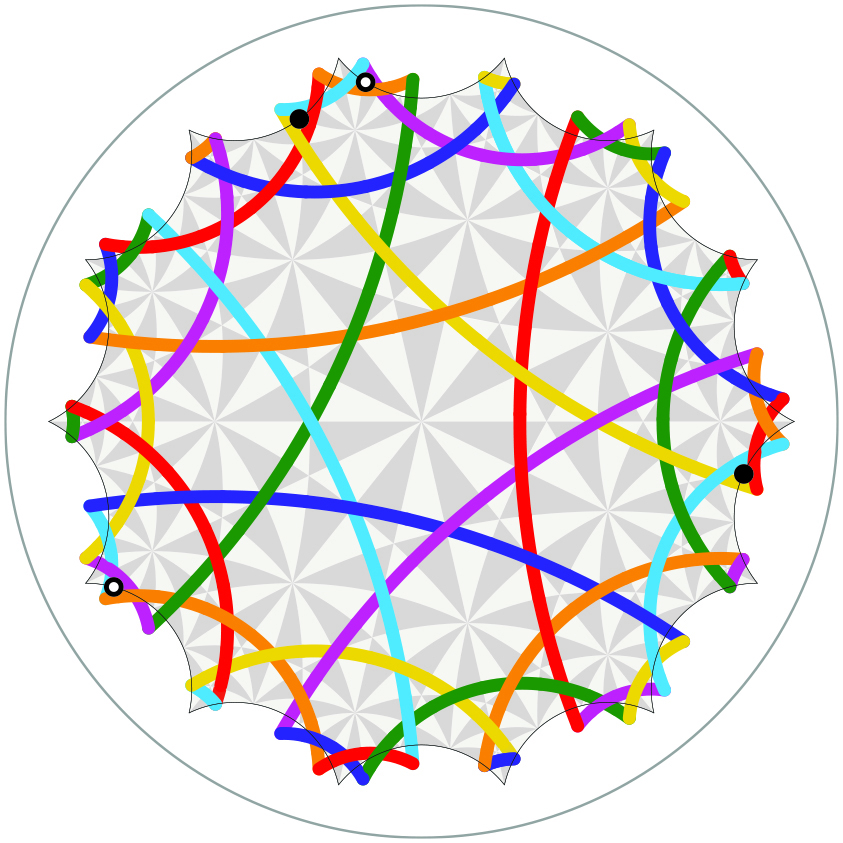}\hspace{.3in}\includegraphics[width=2.4in]{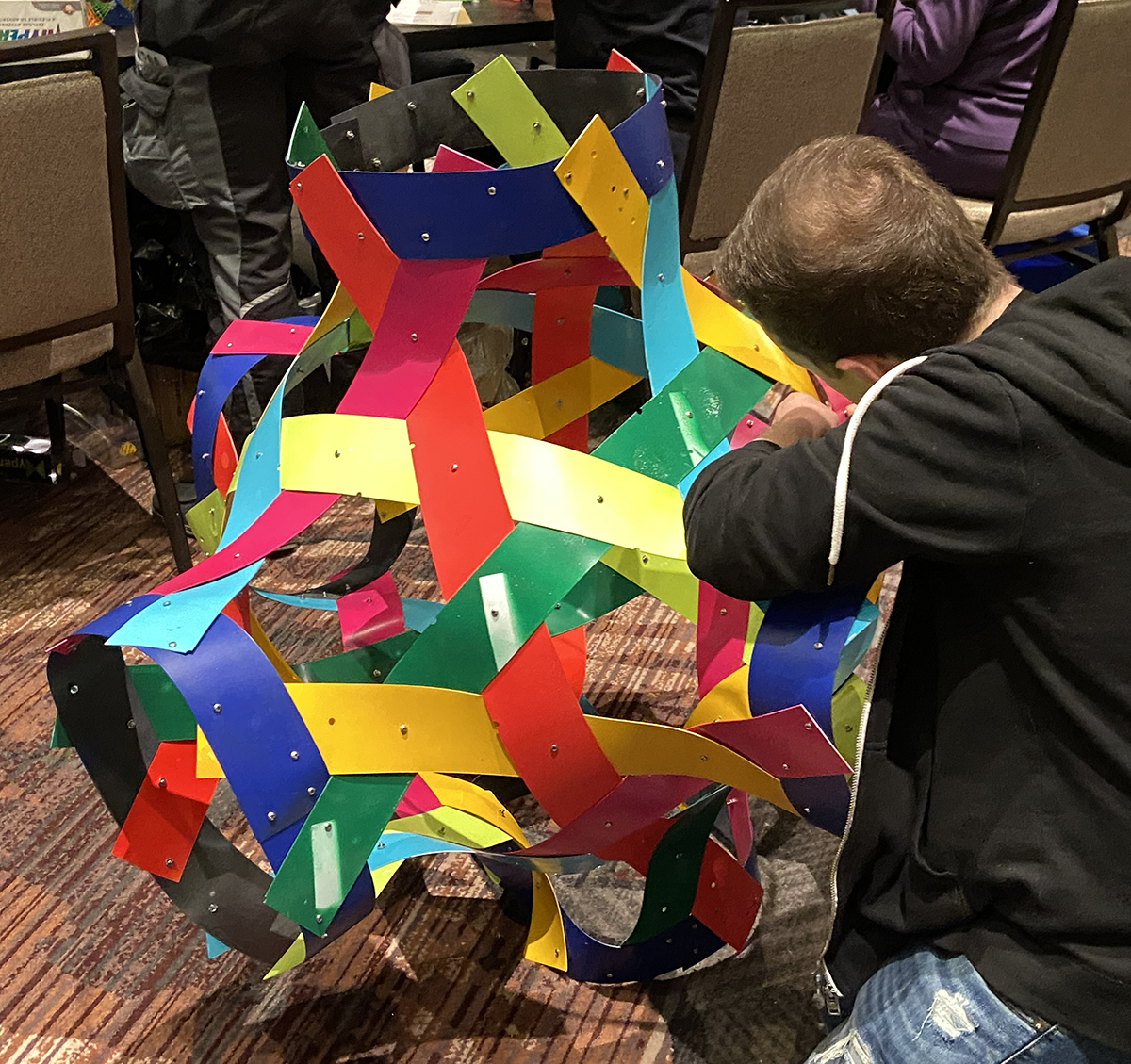}
	\caption{At left a  portion of the hyperbolic plane with underlying *732 symmetry. At right is a model tiled by heptagons and triangles bounded by geodesics, assembled at the Gathering For Gardner, February 24, 2024.  The Klein quartic is formed by identifying like triangles on the boundary of the region at left (two pairs marked), or by identifying opposite blue rings on the model at right.   In both, there are 21 colored geodesics, in seven colors of three each.   Seven triples of colors do not appear, in a Fano plane of the colors, each pair missing from exactly one triple. \label{fig1}}
\end{figure}

\section*{The Klein quartic}

The Klein quartic has a rich history and can be described in many ways, touching many areas of mathematics, as beautifully represented in  {\em The Eightfold Way}, a collection of essays edited by Silvio Levy~\cite{levy}. Models and drawings of it date to its discovery. See~\cite{sabetta, sequin} for a quilted examples and much discussion. 

The original representation of this surface is as a {\em quartic}, the solution set to a fourth-degree algebraic equation in a particular abstract space.\footnote{Specifically, the solutions $[x,y,z]$ in complex projective 2-space  to $x^3y+y^3z+z^3x=0$.} It happens that this solution set has genus-3 and a rich set of 168  orientation-preserving symmetries, the most possible for a surface of this genus. This symmetry group is rather famous; one of its names is PSL(2, 7) and another is PSL(3, 2). We'll call it $Q$ here.

\begin{figure}[h!tbp]
	\centering
	\includegraphics[width=2.2in]{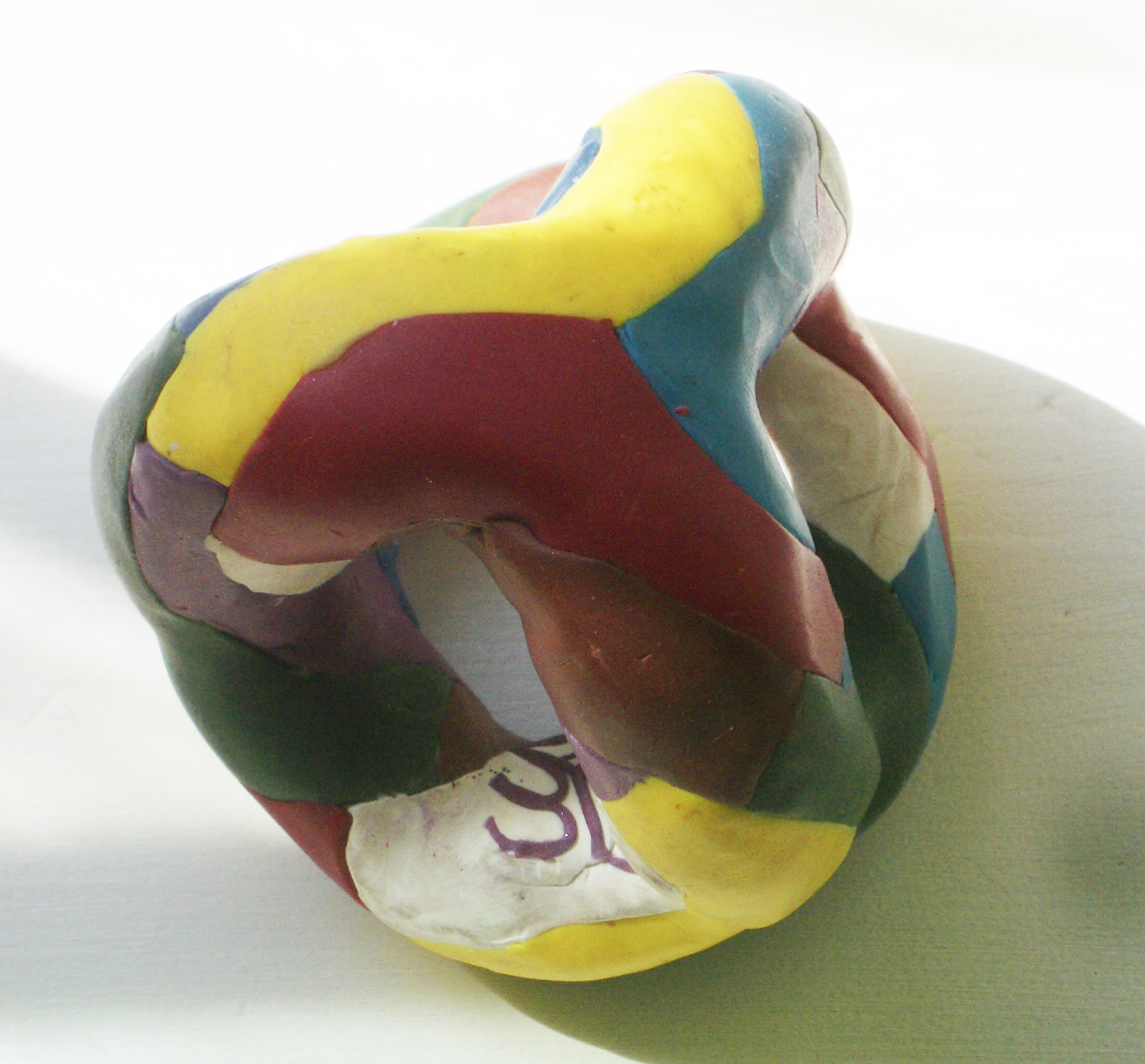}\hspace{.2in}\includegraphics[width=2.8in]{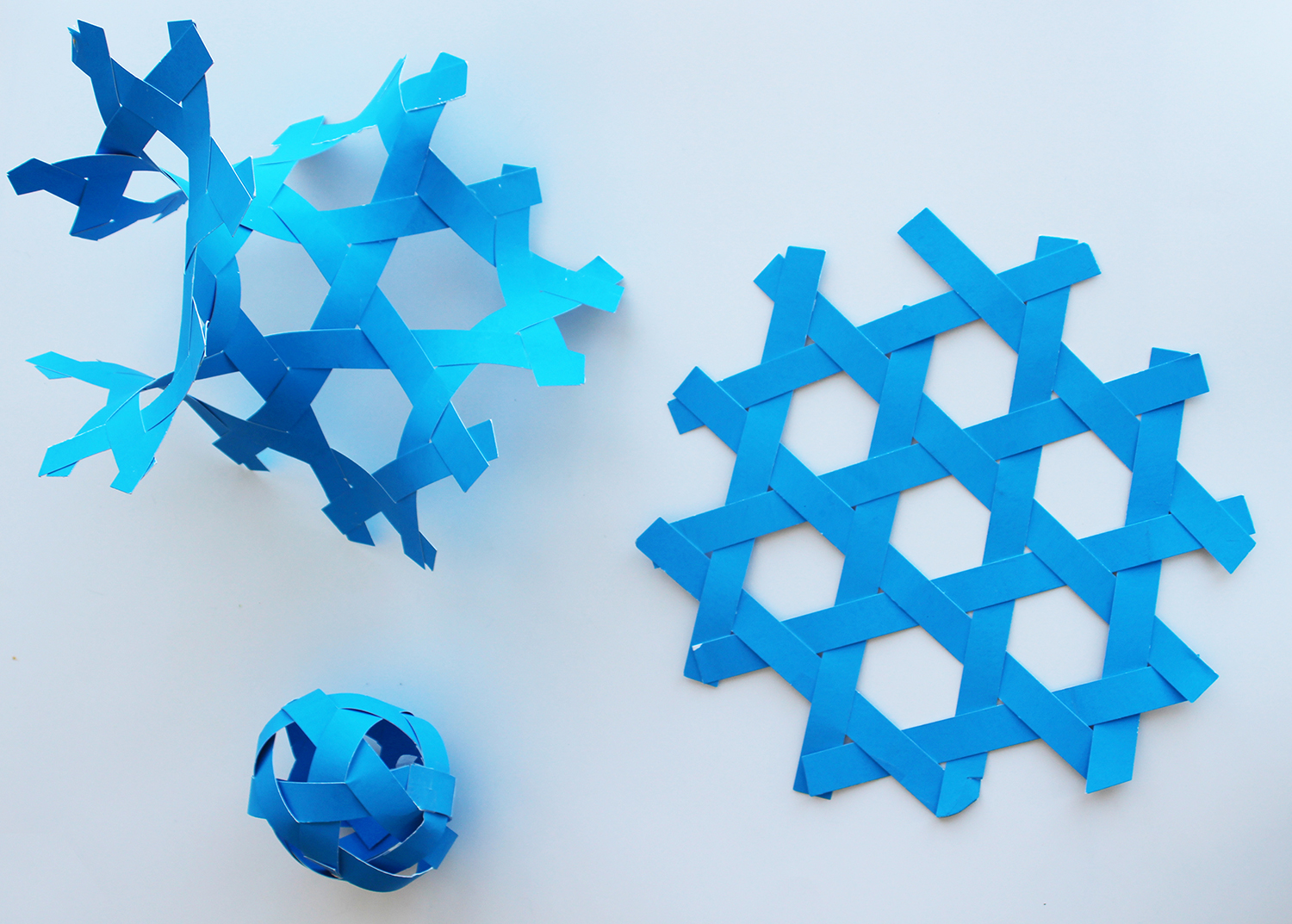}
	\caption{ At left a model of the Klein quartic with tetrahedral symmetry, formed from twenty-four heptagons meeting in threes, in eight colors.
	At right, controlling the surface curvature with the weaving pattern, positive curvature and pentagonal holes at lower right; the negative curvature with heptagonal holes at top; the kogame lattice at right..  \label{fig2}}
\end{figure}

To the  discrete geometer or a low-dimensional topologist, the Klein quartic is most easily understood as a kind of regular polyhedral symmetry, not of the sphere, but as a tiling of a genus-3 surface we'll denote $S$. 
The physical model at right in Figure~\ref{fig2} has much less symmetry, but topologically it has 24 heptagonal facets, meeting  three-to-a-vertex. Abstractly, if we may stretch and deform $S$, then any flag  --- any triple of a  coinciding vertex, edge and face --- may be taken to any other by some topological homeomorphism of the surface, and so this topological tiling is {\em regular}. 
Equally, we can describe the Klein quartic as a regular tiling of a genus-3 surface by 56 triangles meeting in sevens. The two tilings are dual, with the same underlying flags, and $Q$ preserves the tilings and the handedness of the flags.
The term ``Klein quartic'' might refer to either of these tilings, the action of the symmetry group $Q$ on $S$, or any of several other related objects.~\cite{baez,levy}

In his paper~\cite{klein} Klein constructs   the $\{7,3\}$ tiling of the hyperbolic plane, formed by regular  heptagons with $120^{\circ}$ vertex angles, or equally, its dual tiling by equilateral triangles meeting in sevens, both with symmetry group denoted $\ast 732$. He shows how to wrap these tilings onto $S$. The group $Q$ is a quotient of the  orientation preserving symmetry group  $732$. The surface $S$ is the quotient of the hyperbolic plane by the group  $732/Q$. Equivalently, we can assemble the surface by gluing together marked sides of the polygon at left in figure~\ref{fig1} in such a way that colored paths continue and matching triangles on the boundary of the polygon are identified --- two such pairs of triangles are indicated. On the glued-up surface, thick colored lines are  geodesics, three of each color, all of the same length: Eight heptagons  zig-zag across each geodesic before it closes up into a loop.

We are attaching the sides of a 14-sided polygon producing a surface with one fourteen-sided face, seven edges, and (checking carefully) two vertices, for an Euler characteristic of $1-7+2 = -4$. As this surface has no boundary and is orientable, it has genus 3. 
We have thus decorated $S$ with a metric of  constant negative curvature, so that each heptagon is genuinely equilateral and equiangular, all preserved under the action of the symmetry group $Q$.
 
 Physical models must  deform this geometry but we can preserve some symmetry: The  model at right in Figure~\ref{fig2} has tetrahedral symmetry; the group $332$ is a ``subgroup'' of $Q$. (Just as $332$ is a subgroup of $532$, it is not a ``normal'' subgroup of either --- $532\approx A_5$ and $Q$ are ``simple''.)
 S. Matsumoto's fabric model~\cite{sabetta} has more abstract kind of tetrahedral symmetry, in the operations that preserve it. 
 
 We can also represent our genus-3 surface $S$ naturally in our space by cutting it open along three disjoint non-separating curves. We may arrange these six newly cut-open  boundaries  symmetrically  in space, along the usual coordinate axes with matching boundaries opposite one another.  In turn a surface formed from a lattice of these units  is a topological cover of the Klein quartic and separates space into two congruent latticeworks. 
 
 At a glance it would not appear that  twenty-four heptagons could be nicely placed on this surface, but fortunately 
  G. Westendorp shows us how on his website~\cite{westy}   . Once in hand  the example is easier to explain: $Q$ has a subgroup  isomorphic to $432\approx A_4$. This subgroup preserves one of the colors of geodesic. Cutting along these must produce boundaries that are geodesics in the physical model. Consequently they must be flush with the boundary of the unit cell --- they are loops.

As has been beautifully noted before ~\cite{ayres, kagome, agmartin} traditional (Euclidean) basketry and caning  patterns can be adapted to produce surfaces of varying curvature (right, figure~\ref{fig2}). In the ancient Japanese kogame pattern, strips of material are woven in triangular junctures to form hexagonal holes. With pentagonal holes, we have positive curvature, the weaving of a sphere by six bands each at the equator of an icosadodecahedron. With seven, we have negative curvature, and we can realize the pattern shown at left in figure~\ref{fig1} in physical space. (This angle  is actually now close to 58.057 degrees, but there is  enough ``give'' that we can treat it as $60^{\circ}$.)

All of the strips lie along geodesics on the surface of the weaving, at least to the extent that any of that is well-defined. Any bending in the strips must be perpendicular to their normal, and so the strips remain ``straight'' --- and the paths they follow are indeed straight in the underlying hyperbolic tiling.

\section*{How to weave a Klein quartic}

\begin{figure}[h!tbp]
	\centering
	\includegraphics[scale=.75]{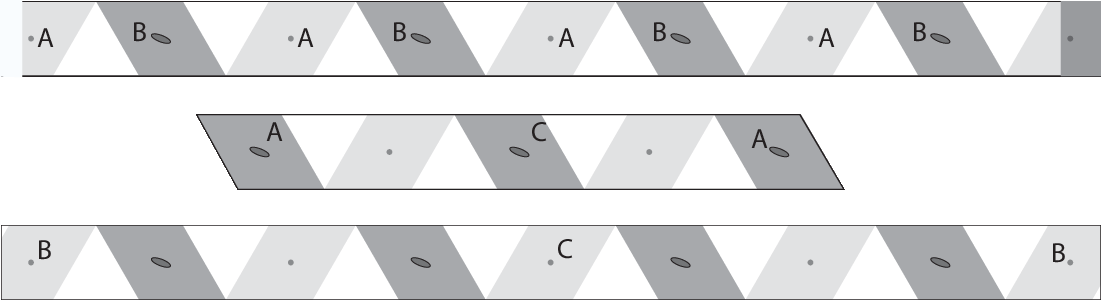}
	\caption{Templates, at top for the ``rings'', middle for the ``short strips'', and bottom for the ``long strips''. 	
	  \label{fig3}}
\end{figure}

With remarkably little time and material you can have a Klein quartic of your very own. 
You will cut a few strips from each of seven colored sheets of papers. Scaling the templates in Figure~\ref{fig3} by 200\% fits US office paper well.

At top in figure~\ref{fig3} is a template for the six ``rings'', all of one color, that will bound the model. Form these by overlapping the left end of the template over the gray region on the right. At middle of the figure is a template for ``short strips" and at bottom is the template for ``long strips''; two short strips and two long strips are needed in each of six other colors, for a total of 18 strips all together.  Dark gray regions of the template will pass under another strip; lighter gray regions will pass over. The  matching letters in the templates are described below. 

\begin{figure}[h!tbp]
	\centering
	\includegraphics[width=\textwidth]{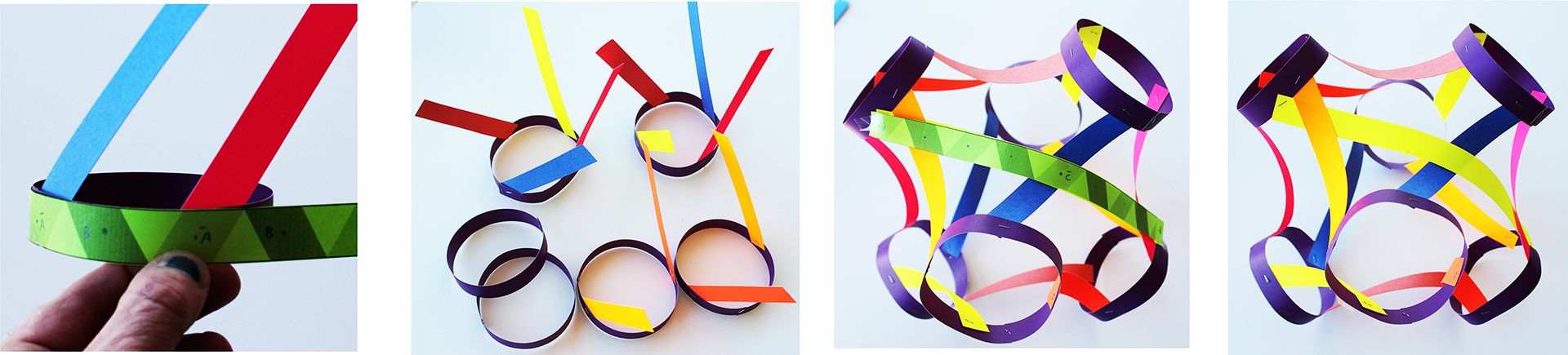}
	\caption{Tips for constructing the model are described in the text.
	  \label{fig4}}
\end{figure}

Begin by attaching four differently colored short strips to the {\em inside} of a ring as shown, using the template to measure their spacing, as at left in figure~\ref{fig4}, matching $A$ on the ring to $A$ on a strip. 

However that ring is colored (say ABCD), attach the same colored strips in opposite order (DCBA) to another ring. Attach each of the two remaining colors to opposite sides of  two of the remaining rings, as shown at middle left in the figure.  

Next, attach the loose ends of the short strips to other rings, following two simple rules: the strips and rings will form an octahedral structure, with the short strips on the edges of this underlying octahedron.  Second, strips of the same color will be on opposite edges. 

We  now weave in the long strips. With a little practice the method becomes natural. Each end $B$ of each long strip will weave over a $B$ on a ring, then pass under a neighboring short strip. A short and long strip meet at their middles $C$, the long strip always over the short one. The rest of the weaving follows. 

We determine the colors of the long strips, at right in figure~\ref{fig4}: Each color of short strip meets four of the rings. The other pair of rings will be connected by a pair of long strips in that color. Because strips of the same color do not cross, this determines which way these strips must travel. (Though the weaving will end up as over-under-over-under, etc, it may not be so during its construction!)

Finally,  the surface $S$ will not be complete until we glue opposite rings to one another, at least in our minds.
Looking at opposite rings, we can see that colors continue one through the gluing, and that opposite rings are attached with a one-quarter turn twist, which we can label with arrows or letters.

\section*{Acknowledgements}

I am endebted to everyone who helped assemble these models over the last few months and  especially  to  G. Westendorp and his regular tessellation website~\cite{westy} --- when the insight came he had the tessellation ready.

    
{\setlength{\baselineskip}{13pt} 
\raggedright				

} 
   
\end{document}